\documentclass[a4paper,12pt]{amsart}
\usepackage{amsfonts}
\usepackage{amsthm}
\usepackage{amscd}

\usepackage{amssymb}
\usepackage[dvips]{graphicx}

\numberwithin{equation}{section}
\newtheorem{theorem}[equation]{Theorem}
\newtheorem{lemma}[equation]{Lemma}
\newtheorem{proposition}[equation]{Proposition}

\newtheorem{corollary}[equation]{Corollary}

\newtheorem{example}[equation]{Example}
\newtheorem{definition}[equation]{Definition}
\newtheorem{remark}[equation]{Remark}

\newcommand{\p}{{\partial}}
\newcommand{\al}{{\alpha}}

\newcommand{\om}{{\omega}}
\newcommand{\eps}{{\varepsilon}}

\newcommand{\la}{{\lambda}}

\newcommand{\rar}{\rightarrow}

\newcommand{\x}{\times}
\newcommand{\ms}{\medskip}

\newcommand{\er}{\mathbb R}

\newcommand{\NI}{{\noindent}}

\newcommand{\bk}{\hfill $\Box$}

\title[Constructions of contact forms on products]
{\bf Constructions of contact forms on products and piecewise
fibered manifolds}

\author{Bogus{\l}aw Hajduk, Rafa\l\ Walczak}

\begin{document}

\maketitle

\begin{abstract}

We study constructions of contact forms on closed manifolds. A
notion of strong symplectic fold structure is defined and we prove
that there is a contact form on  $M \x X$ provided that $M$ admits
such a structure and $X$ is contact. This result is  extended to
fibrations satisfying certain natural conditions. Some examples and
applications are given.

\bigskip
\small{

\noindent {\bf Keywords}: contact form, symplectic fold, open book decomposition

\noindent {\bf AMS classification (2010)}: Primary 53D05,53D10.}

\end{abstract}

\section{Introduction}\label{int}

In this paper we study constructions of contact forms on closed
orientable  manifolds. An intricate question of contact topology is
whether any closed almost contact manifold admits a contact
structure. It is solved positively only in dimensions three and five
\cite{Ma,G1,CPP,EJ}. However, even in low dimensions, this is
usually very non-trivial to construct explicitly a contact form on a
given almost contact manifold.

There are some obvious classes of almost contact manifolds. First of
all, the product of an almost complex manifold $M$ (more generally,
a stably almost complex manifold of even dimension) with a contact
manifold $X$ is almost contact. It is known that $M\x X$ is contact
if $M$ is an orientable surface and $X$ is contact (see \cite{Bo}
for the case of genus $>0$ and \cite{BCS} for $M=S^2$). Our aim and
the principal motivation was to understand the case of $M$ closed
and of arbitrary dimension.

There is a simple case when a contact form exists on $M\x X.$  Let
$(M,\omega )$ be an exact symplectic manifold (i.e., its symplectic
form is exact, $\omega =d\beta)$ and  $\eta $ be a contact form on
$X.$  Then the product form $\beta + \eta$ is contact
\footnote{Formally we should write $p_M^*\beta + p_X^*\eta ,$ where
$p_M, p_X$ are projections, but to simplify the notation we omit
projections. For the same reason wedge signs are omitted in exterior
products of forms.}. Exact symplectic manifolds are necessarily
open, so this cannot be applied directly to closed manifolds.
However,  if $M$ is compatibly decomposed into the sum of exact
symplectic pieces, then there is a formula \cite{GS} which yields a
contact form on $M\x S^1.$ To be a bit more precise, $M$ is assumed
to be a sum of exact symplectic cobordisms which meet at their
convex ends and agree on their common boundaries. Thus $M$ is cut by
a hypersurface and along it the symplectic forms of adjacent pieces
yield a fold. We call such decomposition a strong symplectic fold of
convex type (see Section \ref{prel} for the precise definition and
comments on the formula). Theorem \ref{mainth}  says that the
product of a manifold with strong symplectic fold of convex type
with a contact manifold is contact and it is the base for further
construction and applications. The proof has two main ingredients:
the Giroux - Mohsen \cite{Gi,GM} theorem which states that any
contact form can be deformed to a contact form given by an open book
decomposition, and the heat flow deformation of a confoliation to a
contact form given by Altschuler and Wu \cite{AW}. Open books
together with the Geiges - Stipsicz formula enable us to define a
confoliation and the heat flow applied to it gives a contact form.

Then we extend this theorem in two ways. First, we allow some
bundles over the exact symplectic pieces. In particular, we prove
that there exists a contact form on the total space of any bundle
over a strong symplectic fold with contact fiber if some rather
natural conditions are satisfied. For instance, this holds  for any
bundle over $S^{2n}$ if the structure group preserves the contact
form of the fiber.  Secondly, we show that in $M$ one can allow also
concave folds, i.e., the fold is given by  two concave ends of
symplectic cobordisms (see Section \ref{conc}).

We give a number of examples and applications. They include products
$X\x S^{k_1}\x ... \x S^{k_r}$ provided that $X$ is a  contact manifold
and $k_1+...+k_r$ is even and  $M \times S^{k_1}\x ... \x S^{k_r}$
if $M$ is a strong symplectic fold of contact type and $k_1+...+k_r$
is odd. We show also examples of  homogenous spaces which are
contact but have no invariant contact forms. Moreover, we show that
some surgeries and blowing ups preserve contactness (cf. Proposition
\ref{surgery} and Example \ref{blow}). We describe also a
generalization of the open book construction of contact forms (see
Section \ref{genopen}).

To give a sample  of applications,  consider the following
fillability question. Any  contact form $\lambda$ on $X$ yields the
form $e^t\lambda$ on $X\x [0,1]$ with symplectic exterior
differential (called symplectification of $\lambda ).$ $(X, \lambda
)$ is called fillable, if there exists a symplectic form on a
manifold $W$ with $\p W = X$ equal to  $d(e^t\lambda )$ on a collar
of the boundary. There are obstructions to fillability, in
particular in dimension 3 no overtwisted  form is fillable. However,
there is an interesting and natural weaker question whether the
product form $e^t\lambda + d\phi ,$ where $d\phi$ is the standard
orientation form of $S^1$, extends to a contact form on $W\x S^1.$
In Proposition \ref{exttodisk} we construct such extensions from
$X=S^{2n+1}$ to $D^{2n+2}$ for some forms on $S^{2n+1}.$  If $n=1,$
one can use as $\lambda$ also some overtwisted forms. This  shows
that after multiplying with $S^1$ the obstruction to fillability
disappears, at least for some classes of contact forms. This is a
new proof of a result of Etnyre and Pancholi \cite{EP}. See Section
\ref{furapp} for details.

The constructions of contact forms on bundles over strong symplectic
can be localized. This leads to a class of decompositions into
fibered pieces which are still sufficient to get contactness. In Appendix A we
give a preliminary version of this. We will study such notion
together with its applications in a future paper. Appendix B
contains sample computations in low dimensions performed using
Mathematica.

The authors would like to thank Jonathan Bowden, Diarmuid Crowley
and Andr\'as Stipsicz for interesting comments and pointing out an
incorrect statement in the previous version of this paper.

\section{Preliminaries}\label{prel}

We consider compact, smooth,  orientable manifolds and we want to
find constructions of smooth contact forms on a possibly large class
of manifolds.

Geiges and Stipsicz \cite{GS} gave a formula which yields a contact
form on products $M\times S^1$ for some closed $M$. Let us  describe
their construction in a slightly more general setup. We start with
the definition of a structure which is crucial for our main theorem.

\begin{definition} A strong symplectic fold structure of convex type
on a compact manifold $M$ is a decomposition $M=W_-\cup_NW_+,$ where
$N = W_-\cap W_+$ is a hypersurface in $Int\, M,$ together with
exact symplectic forms $\omega_- = d\gamma_-, \omega_+=d\gamma_+$ on
respectively $W_-,\ W_+,$ such that  the forms satisfy the following
convexity conditions on a tubular neighborhood $N\times [-1,1]$ of
$N$ and at $\p M$:
\begin{enumerate}

\item\label{ssf2} $\gamma_- = e^{t}\lambda$ on $N\times
[-1,0]= N\x [-1,1]\cap W_-$ and $\gamma_+(t)=e^{-t}\lambda$ on
$N\times [0,1]= N\x [-1;1]\cap W_+,$ where $t$ is the  parameter of
$[-1,1]$  and $\lambda$ is a contact form on $N,$

\item the closure of every component of $M-N$ containing a component
of $\p M$ is an exact symplectic cobordism (by (\ref{ssf2}), it is
necessarily convex at the $N-$end), either convex or concave at the
component of $\p M.$
\end{enumerate}

\end{definition}

The hypersurface $N$ is called {\it the fold locus.} The product
$N\x [a;b]$ endowed with the form $d(e^t\la )$ is called the {\it
symplectization} of a contact form $\lambda$ on $M.$ Hence, in the
above definition we assume that on both sides of $N$ we have
symplectizations of $\lambda .$ If the Liouville vector
field of $\omega_{\pm}$ is transverse to $N,$ then one can deform
the symplectic forms to symplectizations.

\NI An obvious example is the double $W\cup (-W),$ where $W$ is a
compact manifold with boundary and $W$ admits an exact symplectic
form satisfying convexity condition  (\ref{ssf2}) at $\partial W$.
Note that a strong symplectic fold does not determine the
orientation, since the orientations given by the symplectic forms on
any two adjacent components of $M-N$ are opposite.

In our terminology we follow Ana da Silva \cite{dS}. She shows that
on any closed stably almost complex manifold there exists a
symplectic fold, i.e. a  2-form which is symplectic everywhere
except for a hypersurface, where the form has fold singularities. A
symplectic fold is globally defined and smooth. It is symplectic
outside a hypersurface and gives opposite orientations on any two
adjacent parts. However, in general the symplectic forms are not
exact and the behavior along the singular hypersurface differs from
what we require for strong symplectic folds. For instance,
symplectic folds do not need give contact forms on the singular
hypersurface.

\begin{theorem} \label{gs} \cite{GS} If $M^{2m}$ admits a
strong symplectic fold of convex type, then $M\times S^1$ is
contact.
\end{theorem}

{\bf Proof.} Let $d\phi$ denote the standard orientation form on
$S^1$ and $p:M\times S^1\rar M$  be the projection. If $\omega_{\pm}
= d\gamma_{\pm},$ then $p^*\gamma_{\pm}+d\phi$ are contact forms
outside $N\times [-1,1]\times S^1.$

Choose smooth  functions $f,g: [-1,1]\rar \er$ such that:

\begin{enumerate} \item $g$ is odd, equal to $1$ near $t=-1,$ equal to $-1$ near
$t=1,$ and it is decreasing  from $-1$ to
$1,$

 \item $f$ is even, positive, equal to $e^{\pm t}$ near
$\pm 1$ and increasing on $[-1,0],$

\item $f'g-g'f>0$ on $[-1;1].$
\end{enumerate} Then the formula

\label{forgs}$$\alpha = f\lambda + g\, d\phi$$ on $[-1,1] \times N
\times S^1 $ yields a contact form on $N\times [-1,1]\times S^1$
(with contact form $\lambda$ on $N$) which extends those defined
above.
  In fact, it is not difficult to calculate:

$$\alpha (d\alpha )^n = nf^{n-1}(f'g-fg')dt  \lambda
 (d\lambda)^{n}  d\theta > 0.$$  \bk

\medskip

Geiges and Stipsicz apply this formula  to show that for every
closed orientable 4-manifold $M$ the product $M\x S^1$
 is contact. They use \cite{B} where it is shown that any closed
 orientable 4-manifold admits a strong symplectic fold of convex
 type.

We want to use the above formula when the circle is replaced
by a general contact manifold $X.$ For this purpose it is necessary
to have a pair of contact forms on $X$ defining opposite orientations
and connected by a path of forms with controlled disruption of
contactness. Moreover, one can see rather easily that the
orientation change should be "one dimensional", for example given by
changing the direction of a vector field transversal to the contact
structure. To construct this we will use contact forms defined in
terms of open book decompositions. So let us recall this
construction.
\medskip

\begin{definition}\label{obd} An open book decomposition of $X$ is given by
\begin{enumerate}
\item a codimension two submanifold $B\subset X$ (called the binding),
\item a tubular neighborhood $U$ of $B$ diffeomorphic to
$B \times D^2,$
\item a fibration $\pi:E= X-B \rightarrow S^1$ with
fiber $P$ (called the page)

\end{enumerate}
such that the monodromy of the fibration $\pi$ is equal to the
identity in $P\cap U$ and   $\pi|U$ can be identified with the
standard projection $B \times (D^2-\{0\}) \rar S^1.$
\end{definition}

\medskip

According to \cite{TW}, one can associate a contact form (which we
will call of {\it open book type}) with any open book decomposition
satisfying the following conditions:

\begin{enumerate}
\item P is exact symplectic, i.e., $P$ has 1-form
$\beta$ such that $d\beta$ is symplectic on $P,$
\item a tubular neighborhood $U$ of $\partial P$ is of convex type,
which means that in a collar $\partial P \x [0,\epsilon)$ we have
$\beta=e^{-t}\nu$ with $\nu$ contact on $\partial P,$
\item the monodromy $f:P\rightarrow P$ of $\pi$ is exact, which means that
$f^*\beta-\beta=d\psi$ for some function $\psi:P\rightarrow
\mathbb{R}.$
\end{enumerate}

Before we write a formula for such form, let us note that the main
theorem of \cite{Gi,GM} says that any contact form is homotopic
(i.e., there exists a deformation through contact forms) to a form
of open book type. It is not unique, but assuming that a contact
form is of open book type does not restrict generality.

If $f:P\rightarrow P$ is the monodromy of $\pi,$ we identify $E$
with the quotient of  $P \times [0,2\pi R]$ for some fixed $R,$ by the
identification $\Phi:(x,0)\sim (f^{-1}(x),2\pi R).$

On $P \x [0,2\pi R]$ we put $\eta_E= \overline{\beta}+d\phi$ with
\begin{eqnarray}\label{beta}\overline{\beta} =
\beta+u(\phi)d\psi\end{eqnarray} for some non-decreasing function
$u:[0;2\pi R]\rightarrow [0;1]$ so that for a small $\varepsilon>0$

\begin{eqnarray}\label{u}u(\phi)=\left\{ \begin{array}{rl} 0 & \textrm{for
$\phi\in[0;\varepsilon)$}
\\ 1 & \textrm{for $\phi\in(2\pi R-\varepsilon;2\pi R].$}
\end{array} \right.
\end{eqnarray}

\NI The form $\overline{\beta}$ descends to $(P \times
[0,2\pi R])\slash \hspace{-1mm} \sim$
since

$\Phi^*(\beta+d\psi+d\phi)=\beta+d\phi$ and $\eta_E$ defines a
smooth form on $E.$ Moreover, if dimension of $P$ is $2n,$ then
$\eta_E  (d\eta_E)^n =
d\phi(d\beta)^n+n\beta(d\beta)^{n-1}u'(\phi) d\phi d\psi.$ As
$d\beta^n>0$ on $P$ and for $R$ big enough the derivative
$|u'(\phi)|$ can be made arbitrary small, $\eta_E$ is contact.

\begin{remark} As far as we know, such "enlarging the circle"
trick  has never been used before in this context. When we tried to
apply the formulae we had been able to find in the literature, then
we needed  an additional assumption, essentially that the fibration
$E\rar S^1$ was trivial. It was rather unexpected that the simple
trick described above enabled us to solve this problem.
\end{remark}

In the sequel we will use a deformation of such form to one having
the opposite orientation of $S^1$ in the fibration $E \rightarrow
S^1.$ For this reason we have to consider the family of forms
$\eta_E=\beta+u d\psi(\phi) +l d\phi$ depending on $l \in
\mathbb{R}.$ Now $\Phi^*(\beta+d\psi+l\cdot d\phi)= \beta+l \cdot
d\phi ,$  so $\eta_E=\overline{\beta}+l\cdot d\phi$ is well-defined
on $(P \times [0,2\pi R])\slash \sim$ for any $l \in \mathbb{R}.$
We have
\begin{eqnarray}\label{eta1}
\eta_E  (d\eta_E)^n = ld\phi\, \left((d\beta)^n-n\beta
(d\beta)^{n-1}u'(\phi) d\psi \right).\end{eqnarray} Note that the formula
implies that our choice of $R$ does not depend on $l$ and we get

\begin{proposition}\label{uphi} If   $R$ is large enough, then
all forms in the family $\eta_E$ are contact for  $l\neq 0.$
\end{proposition}

As the monodromy $f$  is the identity near the boundary $\partial
P,$ the form $\overline{\beta}+l\cdot d\phi$ ($l \in \mathbb{R}$) is
equal to $\nu e^r + l\cdot d\phi$ near the boundary of $B \x D^2$ in
polar coordinates $(r,\phi)$ on $D^2.$ Now we extend
$\overline{\beta}+l \cdot d\phi$ to $B \x D^2$ by the formula

$$\alpha=h_1(r) \nu+l \cdot h_2(r)d\phi ,$$
where

\begin{eqnarray}\label{h1}h_1(r)=\left\{ \begin{array}{rl}
2 & \textrm{for $ r=0$}
\\ e^{1-r} & \textrm{for $r \in [1;R],$}
\end{array} \right.\end{eqnarray}
is strictly decreasing with all derivatives at $0$ vanishing,

\begin{eqnarray}\label{h2}h_2(r)=\left\{ \begin{array}{rl}
r^2 & \textrm{near $r=0$}
\\ 1 & \textrm{for $r \in [1;R]$}
\end{array} \right.\end{eqnarray}
and nondecreasing with $h_1(r) h_2'(r)-h_1'(r) h_2(r)>0$ (see the
drawing below). As another simple calculation shows, the resulting
form is contact on $X.$

\smallskip

\centerline{\includegraphics[width=6cm,height=4cm]{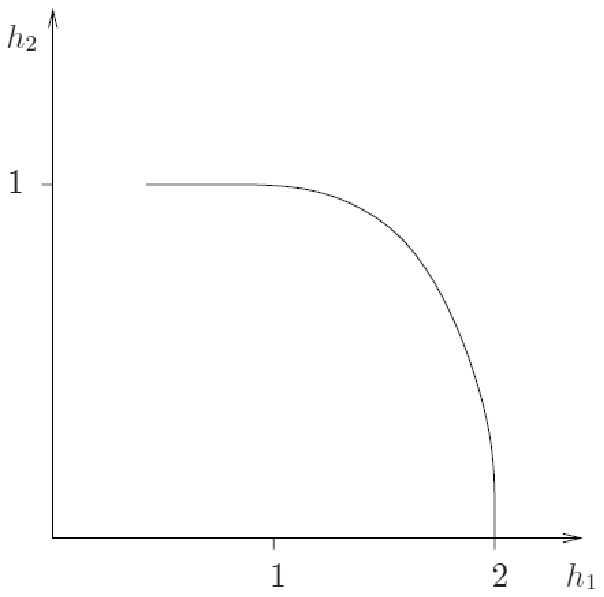}}

If $l=\pm 1,$  for a suitable choice of $u$ and $R$ big enough, both
forms $\eta_E = \overline{\beta} \pm d\phi$ are contact. They
determine opposite orientations and we use this pair of forms
together with the family $\eta_E, l\in \mathbb R,$  in the sequel.

\smallskip

\NI {\bf Notation.} If $\eta$ is one form of such a pair, then by
$\hat\eta$ we denote the other one.

\begin{remark}\label{concordance}  In some proofs in the sequel
we use the following well-known fact: if $\eta_1,\eta_2$ are contact
and homotopic on $X,$ then there is a topologically trivial
symplectic cobordism $M = X \times [0,1]$ between $(X,\eta_1)$ and
$(X,\eta_2).$ Let $M$ be a compact manifold with boundary of contact
type such that the resulting form on $\partial M$ is $\lambda.$ If
we have a homotopy from $\lambda$ to $\lambda',$ we can add a
trivial cobordism to the boundary of $M$ so that we get $\lambda'$
on $\partial M.$ In particular, for a manifold with boundary of
contact type, we can always assume that  we have a contact form of
open book type on $\partial M.$
\end{remark}

Our principal analytic  tool is  the heat flow deformation of a
confoliation \cite{AW}. On a closed manifold $Y^{2m+1}$ consider a
confoliation, i.e. a 1-form $\alpha$ satisfying the inequality
$\alpha \wedge (d\alpha)^{m}\geq 0.$ The points $x \in Y$ where
$\alpha \wedge (d\alpha)^{m}>0$ are called contact (regular), the
other (non-contact) points are called singular and the set of
singular points will be denoted by $\Sigma.$ Altschuler and Wu show
that under some assumptions, the heat flow can deform the
confoliation to a contact form. To describe those assumptions we
choose a Riemannian metric $g$ on $Y$ and consider the form $\tau =
\star(\alpha \wedge (d\alpha)^{m-1}),$ where $\star$ denotes the
Hodge star. Then at every point $x\in Y$ we denote by $\mathcal D
\subset TY_x$ the orthogonal complement of $Null(\tau)_p=\{ V\in
T_pY: \iota_V\tau =0\}.$   At a contact point the subspace $\mathcal
D$ has dimension $2m$ and it is perpendicular to $Null(\tau)_p.$
At a point where rank of $d\alpha$ on $ker\, \alpha$ is
$2m-2,$ the dimension of $\mathcal D$ is 2, and $\dim \mathcal D$ is zero
at points where rank of $d\alpha|ker\, \alpha$ is less than $2m-2.$
A point $x$ is called {\it accessible} if there is a smooth curve
$\sigma :[0,1]\rightarrow Y$ such that $z'(t)\in \mathcal D$ and is
non-zero for all $t \in [0,1],$ $z(0)=x$ and $z(1)$ is a contact
point. Thus we see that in the case when the rank of $d\alpha|ker\,
\alpha$ is less than $2m-2$ no singular point is accessible. Since
we have to reduce the general case to that of corank at most 3, this
is one of the main difficulties of our construction.

In the sequel we will use the following theorem.

\begin{theorem}\cite{AW}\label{conf} Suppose that $Y$ is a closed manifold with
a confoliation $\alpha.$ If every non-contact point of $Y$ is
accessible, then $Y$ supports a contact form  $C^{\infty}$-close
to $\alpha$.
\end{theorem}

\section{Main theorem}\label{mainsec}

Our main theorem is the following.

\begin{theorem}\label{mainth} If $(X^{2m+1},\alpha )$ is a
closed contact manifold and $M^{2n}$ admits a strong symplectic fold
of convex type, then $X\x M$ is contact.
\end{theorem}

{\bf Proof.} Consider the decomposition $M=W_1 \cup (N \times[-1;1])
\cup W_2$ and the forms $\omega_+, \omega_-, \lambda$ given by the
strong symplectic fold on $M.$ Here $N=W_+ \cap W_-$ and $N\x
[-1,1]$ is a tubular neighborhood of $N$ with $N\x [-1,0]\subset
W_-, N\x [0,1]\subset W_+,$ $W_1=W_--N\x (-1,0], W_2= W_+-N\x
[0,1),$ $\omega_{\pm} = d\gamma_{\pm}.$

We can assume that the contact form $\alpha$ is of open book type
with  $P,B$ denoting the page and the binding. We use the notation
introduced in Section \ref{prel}: $\beta, \psi$ (see the
paragraph right after Definition \ref{obd}), function $u$ (formula
(\ref{u})), and $h_1, h_2$ (formulae (\ref{h1}),(\ref{h2})),
$\overline{\beta}=\beta+u(\phi)d\psi$ (formula (\ref{beta})).
Recall that $d\beta$ is exact symplectic
on $P$ and a tubular neighborhood of $\partial P$ is of convex type.

We define a 1-form $\tilde\eta$ on $X\x M$ by separate formulae on
$X\x (W_1\cup W_2),\ (X-B\x D^2)\x N\x [-1,1],$ and $B\x D^2\x N\x
[-1,1].$

On $X\x W_{\pm}$ we take $\tilde\eta = \alpha + \gamma_{\pm}.$ By
the discussion in Section \ref{prel}, for every $l\in \mathbb{R}$
the form $\overline{\beta}+l\ d\phi$ is well-defined on $X-B \times
D^2,$ where $D^2$ denotes disk of radius $R.$ This constant $R$ will
be chosen later. Hence $\tilde{\eta}
=\beta+u(\phi)d\psi+g(t)d\phi+f(t)\lambda$ with $f,g$ given in
Theorem \ref{gs}, $t\in [-1;1]$ is well-defined.

Finally, let

\begin{eqnarray}\label{defeta}
\tilde\eta=h_1(r)\nu+f(t)\lambda+h_2(r) g(t)
d\phi \end{eqnarray}

\NI on $B\times D^2\times N \times [-1,1].$

Here $(r,\phi)$ are polar coordinates on the disk $D^2$ and $f,g$
are functions defined in Theorem \ref{gs}.

\begin{lemma}\label{extent} The form $\tilde\eta$ is well-defined
and smooth on $X\times M,$ contact in the complement of
 $(B \times \{0\}) \x (N \times \{0\}) \subset X \x M.$
\end{lemma}

{\bf Proof.} The form $\tilde\eta = \alpha + \gamma_{\pm}$ on $X\x
W_{\pm}$ is clearly contact. By the assumptions and the choices we
made, the partial forms agree smoothly, so we get a globally defined
smooth form.

In  $(X-B \times D^2) \x N \x [-1;1]$ we have
$$\tilde\eta (d\tilde\eta)^{m+n}=
\left(m+n \right)\left(\!\begin{array}{c}m+n-1 \\ m \end{array}\!\right)f^{n-1}(f'g-g'f)(d\beta)^m
d\phi \lambda (d\lambda)^{n-1} dt+$$
$$+\left(m+n \right)\left(\!\begin{array}{c}m+n-1 \\ m \end{array}\!\right)f^{n-1}(d\beta)^{m}
(d\lambda)^{n-1}\left(g'(t)
dt d\phi+f'(t)dt \lambda\right)u(\phi)d\psi+u'(\phi)\kappa$$ with
$2m+2n+1$-form $$\kappa=2\left(\!\begin{array}{c}m+n \\ 2 \end{array}\!\right)
\left(\!\begin{array}{c}m+n-2 \\ m-1 \end{array}\!\right)f^{n-1}f'\beta(d\beta)^{m-1}
(d\lambda)^{n-1}d\phi d\psi dt\lambda.$$

 As $d\beta^m d\psi$ is a
form on $P,$ hence it vanishes  for dimensional reasons and  the
middle term  of the formula for $\tilde\eta (d\tilde\eta)^{m+n}$ is zero. Furthermore, for $R$ big enough
$|u'(\phi)\kappa|$ can be made arbitrarily small, because $\kappa$
does not depend on $R.$ It follows that
$\tilde\eta(d\tilde\eta)^{m+n}>0,$ hence our formula defines a
contact form on this part.

\NI It remains to examine $\tilde\eta$ on $B \x D^2 \times N \times
[-1;1].$ Direct computations give
$$\tilde\eta(d\tilde \eta)^{m+n}=
c_1\left(f'g(h_1h_2'-h_1'h_2)+fg'h_1'h_2\right) \nu  d\nu^{m-1}
 dt \lambda  d\lambda^{n-1} dr d\phi
,$$ where $c_1$ is a  positive constant. Since $h_1h_2'-h_2h_1'>0,$
$f'g\geq 0$ by definition given in Theorem \ref{gs}, and
consequently $fg'h_1'h_2 \geq 0,$ we see that
$\tilde\eta(d\tilde \eta)^{m+n}\geq 0$ and it vanishes if and
only if $f'g=0$ and $fg'h_1'h_2=0.$  The equality $f'g=0$ implies
$t=0.$ Furthermore, for $t=0$ we have $fg'>0.$  To complete the
proof notice that our assumptions on $h_1,h_2$ (i.e.
$h_1'(r)=0\Leftrightarrow r=0, h_2(r)=0 \Leftrightarrow r=0$) yield $h_1'h_2 =0
\Leftrightarrow r=0.$  \bk

We want to apply Theorem \ref{conf}, so we need the accessibility
condition to be satisfied. We know (\cite{AW}) that the necessary
condition for accessability is that $rank \ d\widetilde{\eta}\mid
\ker \widetilde{\eta}\geq 2(m+n-1)$
  Unfortunately,  on $\Sigma = B \x \{0\} \x N \x \{0\}$  we have
$rank \ d\widetilde{\eta}\mid \ker \widetilde{\eta} <2(m+n)-2$ since
 $d\widetilde{\eta}|T(X \x M)|_{\Sigma}=2d\nu+d\lambda$ and
$\widetilde{\eta}|T(X \x M)|_{\Sigma}=2\nu+\lambda.$ Thus the
singular points are not accessible. In order to remedy this we
change the confoliation form making it asymmetric with respect to
the decomposition  $W_1 \cup (N \times[-1;1]) \cup W_2.$ Roughly
speaking, we impose in this way some more transversality along the
singular set. Define the form $\eta$ on $X \times M$ by the formula

\begin{eqnarray}\label{main} \eta=\left\{ \begin{array}{rl}
e^{-1}(\overline{\beta}+d\phi+\gamma_{-}) & \textrm{on $B \x D^2 \x  W_1$}
\\ k(t)\left(h_1(r)\nu+f(t)\lambda+h_2(r) g(t) d\phi\right) & \textrm{on $B \x D^2 \x
N \times [-1;1]$} \\ e(\overline{\beta}-d\phi+\gamma_{+}) & \textrm{on $B \x
D^2 \x W_2.$}
\end{array} \right.\end{eqnarray}

\NI In formula above $k:[-1-\varepsilon,1+\varepsilon]\rightarrow [e^{-1};e]$ is a
smooth, positive, non-decreasing function satisfying

$$k(t)=\left\{ \begin{array}{rl}
e^{-1} & \textrm{on $(-1-\varepsilon;-1]$}
\\ e^t & \textrm{on $[-1+\varepsilon;1-\varepsilon]$}
\\ e & \textrm{on $[1;1+\varepsilon)$}
\end{array} \right.$$
\NI with $\varepsilon$ small enough.

\NI Because $k>0$ and $\tilde{\eta}$ is contact on the complement of $\Sigma,$  hence
$\eta$ is also contact on $X\x M - \Sigma.$  By continuity, $\eta (d\eta)^{m+n}\geq 0$ on $\Sigma.$
Therefore we get again a smooth confoliation with the same critical
set $\Sigma = B\x\{ 0\} \x N\x \{0\} .$

 To apply  \cite{AW} we choose a Riemannian metric
$\langle\cdot,\cdot\rangle$ on $X \times M$ such that near $\Sigma$
submanifolds $N,I=[-1,1],B,D^2$ are pairwise orthogonal. We will
check that $\eta$ satisfies the assumption of Theorem \ref{conf}.
Consider  $\tau = \star\left(\eta (d\eta)^{m+n-1}\right)$ and
$\mathcal D = Null(\tau )^{\perp}.$ We will show that for every
point $(b,v) \in B\times \{0\} \x N \times \{0\} ,$ the radial path
$z(r)=(b,(r,\phi),0,v) \subset B \times D^2 \times N \times I$ (with
$z'(r)=\frac{\partial}{\partial r} \in TD^2$ for $r\in (0;R]$ and
any fixed $\phi \in [0,2\pi)$) satisfies $z'(t)\in \mathcal D,$
hence every $x \in \Sigma$ is accessible from a contact point. The
proof is divided into two parts. We check first that we have
$\mathcal{D}=TD^2$ on $\Sigma$ and then that $z'(r) \in \mathcal{D}$
for $r\in (0;R].$

\begin{lemma}\label{D2} Under the assumptions above, $\mathcal{D}=TD^2$ on
$\Sigma.$
\end{lemma}

 {\bf Proof.} By formula
(\ref{main}), $\eta=e^t\widetilde{\eta},$ hence
\begin{eqnarray}\label{der} d\eta=e^tdt\widetilde{\eta}+e^td\widetilde{\eta}=
\end{eqnarray} $$=e^tdt\left( h_1(r)\nu+f(t)\lambda+h_2(r) g(t) d\phi \right) +$$
$$+e^t \left( h_1'(r)dr\nu+h_1(r)d\nu+f'(t)dt\lambda
+f(t)d\lambda+h_2'(r)g(t)drd\phi+h_2(r)g'(t) dtd\phi \right) .$$

\smallskip

Substituting $t=r=0$  gives that $\eta|T(X \x M)|_{\Sigma}=
2\nu+\lambda$ and $d\eta|T(X \x M)|_{\Sigma}=2d\nu+d\lambda +
dt(2\nu+\lambda)$ with $\Sigma=B \x \{0\} \x N \x \{0\}.$ As
$d\nu^m=0,d\lambda^n=0$ on $\Sigma ,$ we easily calculate:
$$\eta (d\eta)^{m+n-1}=\eta (m+n-1)
(2d\nu+d\lambda)^{m+n-2} dt
(2\nu+\lambda)=$$
$$=\eta(m+n-1)2^{m-1}(d\nu)^{m-1}
(d\lambda)^{n-1} dt (2\nu+\lambda)=$$
$$=C \nu (d\nu)^{m-1} \lambda
(d \lambda)^{n-1} dt=Cdvol_{B} dvol_{N} dt$$ for
some positive constant $C.$ Thus $\star\left(\eta
(d\eta)^{m+n-1}\right)=\pm C dvol_{D^2}$ and $\mathcal D = TD^2.$\bk

\medskip

The last lemma implies that $z'(0)\in \mathcal{D}.$ However, it is
not clear yet if  $z'(r)\in \mathcal{D}$ beyond $\Sigma .$ So now we
determine $\mathcal{D}$ for $r>0.$ The proof is an elementary but
long computation, hence we skip some parts of it.

Since $\eta$ is contact on $X \x M - \Sigma ,$ by \cite{AW} we have
that 2-form $\tau = \star\left(\eta (d\eta)^{m+n-1}\right)$ is of maximal rank $(=2(m+n)),$ and its nullity bundle
$Null(\tau )$ is 1-dimensional. Thus $\mathcal{D}=Null(\tau )^{\perp
}$ is $2(m+n)$-dimensional. For the remaining part of the proof, it
is enough to check that $Null(\tau )$ is perpendicular to
$\frac{\partial}{\partial r}$ on $B \times (D^2- \{0\})\times N
\times \{0\}.$ By our choice of metric, $\frac{\partial}{\partial
r}$ is perpendicular to $T_r = B \times S^1_r \times N \times I$
(with $S^1_r=\{p \in D^2: \ |p|=r\}, r>0$). Therefore once we show
that for $t=0$ the bundle $Null(\tau )$ is tangent to $T_r$ or,
equivalently, that on $T_r$ the form  $\tau$ is degenerate (i.e., of
rank $< 2(n+m)= dim\, T_r),$  the proof of Theorem \ref{mainth} is
completed.

\medskip

As in Lemma \ref{D2}, substituting $t=0$ in formula (\ref{der}) gives
$\widetilde{\eta}|T(X \x M)|_{S}=h_1(r)\nu+2\lambda$ and
$d\widetilde{\eta}|T(X \x M)|_{S}=h_1' dr
\nu+ h_1d\nu+d\lambda-h_2 dt d\phi$ on $S= B \x D^2 \x N \x \{0\}.$ We
obviously have
$(d\eta)^{m+n-1}=(dt\widetilde{\eta}+d\widetilde{\eta})^{m+n-1}=(d\tilde
\eta)^{m+n-1}+(m+n-1)(d\tilde
\eta)^{m+n-2}dt\tilde \eta$ on $S.$ Further, as $d\nu^m=0,d\lambda^n=0$ we get
$$(d\tilde \eta)^{m+n-1}=
\left(\!\begin{array}{c}m+n-1 \\ n-1 \end{array}\!\right)
(d\lambda)^{n-1} (h_1' dr \nu+ h_1d\nu-h_2 dt d\phi)^m+$$
$$+\left(\!\begin{array}{c}m+n-1 \\ n-2 \end{array}\!\right)
(d\lambda)^{n-2}(h_1' dr \nu+ h_1d\nu-h_2 dt
d\phi)^{m+1}=$$
$$=(d\lambda)^{n-1}\left((d\nu)^{m-1}(D_1 dr\nu+D_2
dtd\phi)+D_3(d\nu)^{m-2}dr\nu dt
d\phi\right)+$$ $$+D_4(d\lambda)^{n-2}(d\nu)^{m-2}dr\nu dt d\phi$$
\smallskip
for some functions $D_i \ (i\in\{1,2,3,4\})$ of variable $r.$
In a similar manner we calculate $dt\tilde \eta(d\tilde \eta)^{m+n-2}:$

$$dt\tilde \eta(d\tilde \eta)^{m+n-2} = dt\tilde
\eta(h_1' dr \nu+ h_1d\nu+2d\lambda-h_2 dt
d\phi)^{m+n-2}=$$ $$=dt\tilde \eta(h_1' dr \nu+
h_1d\nu+2d\lambda)^{m+n-2}=$$ $$ =dt\tilde
\eta\left(\!\begin{array}{c}m+n-2 \\ n-1
\end{array}\!\right) \left((h_1d\nu)^{m-1}
(2d\lambda)^{n-1}+(m+n-2)(h_1d\nu+2d\lambda)^{m+n-3}h_1' dr
\nu\right).$$

\smallskip

\NI After arduous, but elementary computation we get that $$\eta
 (d\eta)^{m+n-1}=
C_1\nu(d\nu)^{m-1}dr\lambda(d\lambda)^{n-1}+C_2\nu(d\nu)^{m-1}
d\phi(d\lambda)^{n-1}dt+$$
$$+C_3(d\nu)^{m-1}d\phi\lambda(d\lambda)^{n-1}dt+
C_4\nu(d\nu)^{m-2}dr
d\phi\lambda(d\lambda)^{n-1}dt$$
$$+C_5\nu(d\nu)^{m-1}dr
d\phi\lambda(d\lambda)^{n-2}dt+C_6\nu(d\nu)^{m-1}
\lambda(d\lambda)^{n-1}dt$$

\smallskip

\NI for some functions $C_i, \ i=1,\ldots,6$ of variable $r.$
Furthermore,
 $\check{\nu}=\star\left(\nu(d\nu)^{m-2}\right)$ in $B$ and
 $\check{\lambda}=\star\left(\lambda(d\lambda)^{n-2}\right)$ in $N$
 both have maximal ranks  equal to respectively
 $2m-2$ and $2n-2.$ If we additionally set
 $\nu_1=\star\left((d\nu)^{m-2}\right)$ in $B$ and
 $\lambda_1=\star\left((d\lambda)^{m-2}\right)$
 in $N,$ then
$$\tau=\star\left(\eta (d\eta)^{m+n-1}\right)=E_1dtd\phi+E_2 \lambda_1 dr+E_3
 dr \nu_1+E_4\check{\nu}+E_5 \check{\lambda}+E_6 dr d\phi$$
 again for some functions $E_i, \ i=1,\ldots,6$
 of variable $r.$ The pullback of $\tau$ to $T_r=B
\times S^1_r \times N \times I$ via the inclusion $j:T_r
\hookrightarrow M \times X$ yields
$$j^*\tau=j^*\left(\star\left(\eta (d\eta)^{m+n-1}\right)\right)=
E_1dtd\phi+E_4\check{\nu}+E_5 \check{\lambda}.$$ The rank of this
form is equal to $2(m-1)+2(n-1)+2=2(m+n)-2<2(m+n)=\dim T_r,$ hence
$\tau\upharpoonright T_r$ is degenerate on $T_r.$ As we said
earlier, this implies that $Null(\tau )$ is tangent to $T_r,$ hence
$\frac{\partial}{\partial r} \in Null(\tau )^\perp=\mathcal{D}$ for
$r>0.$ This completes the proof.

\bk

\section{Contact forms on bundles}\label{confib}

In this section we discuss constructions of contact forms on bundles
of two types:

\begin{enumerate}
\item \label{exbun} exact bundles:  bundles over a contact base
with exact symplectic fiber and structure group of exact
symplectomorphisms;

\item \label{conbun}contact bundles:  bundles over a strong
symplectic fold of exact type with contact fiber and structure group
of strict contactomorphisms.

\end{enumerate}

In both cases our results require some further properties of the
bundles. For exact  bundles of type \ref{exbun} we will need  the
following property. Let $E\rar B$ be a smooth bundle with fiber $F$
and the structure group $G\subset Diff(F).$ We say that it is  {\it
defined on a hypersurface} $H\subset B$ if its restriction to $B-H$
is trivial and there is a map $a:H\rar G$ such that the map $A:H\x
F\rar H\x F:(x,v)\mapsto (x,a(x)v)$ is  smooth and the bundle is
obtained by gluing the product pieces along $H$ with $A.$ The
definition applies also in the case when $B-H$ is connected. If $B$
is the circle, then as the hypersurface one can take a single point.

Given an exact symplectic manifold $(M, \omega =d\beta),$ denote by
$Ex(M,\beta )$ the group of exact symplectomorphisms and by $Ex(M,\p
M, \beta )$ the group of exact symplectomorphisms equal to the
identity near the boundary.

\begin{proposition}\label{exactbundles} If $\pi : E\rar B$ is a bundle with
compact contact base $(B,\mu ),$  compact exact symplectic fiber
$(F, \omega = d\beta ),$ the structure group contained in the group
of exact symplectomorphisms $Ex(F,\beta )$ and defined on a
hypersurface $H\subset Int\, B,$  then $E$ admits a contact form. If
the structure group is contained in $Ex(F,\p F,\beta ),$ then the
contact form can be chosen to be equal to the product form $R\mu + \beta$
on a collar of $B\x \p F,$ where $R$ is a large enough constant.
\end{proposition}

Proof. Let $A: H\x F\rar H\x F$ be the gluing diffeomorphism. By
assumptions, for any $x\in H$ we have $(A|\{ x\}\x F)^*\beta = \beta
+ d\psi_x,$ where $\psi_x \in C^{\infty}(F).$ Actually, there
exists a smooth function $\tilde\psi$ on $H\x F$ such that this
equality holds with $\psi_x = \widetilde{\psi}(x,\cdot ).$
Consider a tubular neighborhood $U\cong H\x [-1,1]$ of $H.$ For any
positive constant $R$ the form $R\mu + \beta$ is contact on
$\pi^{-1}(B-U)\cong (B-U)\x F.$  On $U$ consider the form

$$\eta = R \mu + \beta + ud\widetilde{\psi}, $$
where $u:[-1,0]\rar [0,1]$ is given by formula (\ref{u}). If dimension
of $F$ is $2m$ and dimension of $B$ is $2n,$ then
$$\eta d\eta^{n+m-1}= C_1R^m\mu d\mu^{m-1}d\beta^n+
C_2R^{m-1} u'\mu d\mu^{m-2}d\beta^n dt d\widetilde{\psi}\,+$$
$$C_3R^{m-2} u' d\mu^{m-1}\beta d\beta^{n-1}dt d\widetilde{\psi}+
C_4R^{m-1} d\mu^{m-1}d\beta^n d\widetilde{\psi},$$ where $C_i,
i=1,2,3,4$ are constants depending only on $m,n.$ For $R$ large
enough the first term dominates the whole sum and consequently
$\eta$ is contact. By construction, the forms on $\pi^{-1}U$ and on
$\pi^{-1}(B-U)$ agree near $H\x \{ \pm1\} \x F,$ hence we obtain
a smooth contact form on $E.$ \bk

\ms

In case of contact bundles consider first bundles over an exact
symplectic manifold.

\begin{proposition}\label{ctctoverex}  Let $(W, \omega_0=d\beta )$ be a compact exact
symplectic manifold,
 $\pi :E\rar W$ a bundle over W with a closed  contact
 fiber $(X,\eta_0 ).$
If the structure group of the bundle is contained in the group $Cont
(X, \eta_0)$ of diffeomorphisms preserving the contact form $\eta_0$
(strict contactomorphisms), then E admits a contact form. If the
bundle is trivial over $\p X$ and $\beta = e^t\lambda$ in a collar
$U$ of $\p W,$ where $\lambda$ is a contact form on $\p W,$ then one
may require the form to be the product form $Ke^t\lambda + \eta_0$
in a collar of $U\x X,$ where  $K$ is a large enough positive
constant.
\end{proposition}

{\bf Proof.} We will use the symplectization of the fiber and the
well-known Thurston construction of symplectic forms on bundles. Let
$\{U_s\}_{s \in S}$ be an open cover of $W$ with local
trivializations $\Psi_s: {\pi}^{-1}(U_s) \cong U_s \times X.$ If
$\{f_s\}_{s \in S}$ is the partition of unity subordinated to
$\{U_s\}_{s \in S},$ then we define a symplectic form
$\omega=d(K\pi^*\beta +e^t(\sum_{s \in S} f_s\Psi_s^*\eta_0))$ on $E
\times  [-\varepsilon, \varepsilon ]$ for some $K$ big enough and
$\varepsilon>0.$ Let $\mathcal R$ be the Reeb vector field of
$\eta_0.$ Its interior products with $\eta_0, d\eta_0$ are
$\iota_{\mathcal R}\eta\equiv 1, \iota_{\mathcal R}d\eta_0 \equiv 0.$ Since $\eta_0$ is
preserved by the structure group of the bundle, there is a
horizontal  vector field $\tilde{\mathcal R}$ on $E$ such that its
pushforward by $\Psi_s$ is equal to $\mathcal R$ for any $s\in S.$ This
implies that $\tilde{\mathcal R}$ is the Reeb field of
$\Psi_s^*\eta_0|\pi^{-1}(w) $ for any $s$ and $w\in W.$ Thus, if
$\eta = \sum_{s \in S} f_s\Psi_s^*\eta_0, $ then we have
$\iota_{\tilde{\mathcal R}}d\eta \equiv 1, \iota_{\tilde{\mathcal
R}}\eta \equiv 0.$ Therefore for the Liouville vector field $L$ of
$\om$ we have $\iota_L\om=K\pi^*\beta + e^t\eta.$ If we additionally
apply $\iota_{\tilde R}$ to the last equation, we get
$-\iota_L\iota_{\tilde R}\om = -\iota_L\iota_{\tilde
R}(e^tdt\eta+e^td\eta)=e^t\iota_Ldt=e^t.$ This implies that $L$ is
transversal to $E,$ hence $E \cong E \times \{0\} \subset E \times
[-\varepsilon ,\varepsilon ]$ is contact. The additional convexity
property of $\p X$ follows from the fact that one can take $U$ as a
trivialization chart (i.e. $U \in \{U_s\}_{s \in S}$). \bk

\ms

Consider now contact bundles over a strong symplectic fold. We will
prove a generalization of Theorem \ref{mainth} in this case.

\ms

\NI Let $(X,\eta )$ be a closed contact manifold and let $(W_{\pm},
N, \lambda_{\pm})$ be a strong symplectic fold of convex type on
$M.$ Consider a bundle $E\rar M$ with fiber $X$ and let $E_{\pm}\rar
W_{\pm}$ denote its restrictions to $W_{\pm}.$  We assume that the
bundle is trivial over the fold locus $N=\partial W_- \cap\partial
W_+,$  $E_-$ is contact with respect to $(X,\eta ),$ $E_+$ is
contact with respect to $(X,\eta' ).$ We also assume that there exists a
contact form $\eta_0$ of open book type on $X$ such that  $\eta$ is
homotopic to $\eta_0$ and $\eta'$ to $\hat\eta_0.$

\begin{theorem}\label{mainthfib} If $E$ is the total space of contact
fibration  over a strong
symplectic fold of exact type and  satisfies the above assumptions,
then there exists a contact form on $E.$
\end{theorem}

For $E_{\pm}$ we apply Proposition \ref{ctctoverex}. Over the collar
$N\x [-1,1]$ the bundles are product, thus the arguments used in the
product case work. To be more precise, we start from the contact
forms on $E_+$ given by the contactness of those bundles. Since the
bundles are trivial over $N\x [-1,1],$ we can use the homotopies
$\eta \sim \eta_0, \eta'\sim \hat\eta_0$ to get the form $\beta
+\eta_0$ over $N\x [-\eps, 0]$ and $\beta +\hat\eta_0$ over $N\x
[0,\eps ]$ for some $\eps
>0.$ Having established this, we can apply the same arguments which
were used to prove Theorem \ref{mainth}.

\begin{corollary}\label{ctctoversph} If one of the fibrations
$E_{\pm}\rar W_{\pm}$ is trivial, then $E$ admits a contact form.
In particular, this holds for any contact fibration over a sphere
with its standard strong symplectic fold  $S^{2n}=D^{2n}_-\cup
D^{2n}_+.$
\end{corollary}

\section{Some applications}\label{app}

Results of the previous sections give a constructive way to show
that some manifolds  are contact.  We present now a series of
examples. First of all we discuss the problem of existence
of strong symplectic folds on manifolds which in general seems to be
a difficult question.

 Let us recall that $W$ is the {\it trace} of
a (single) {\it surgery} of index $k+1$ on $M^{2n+1}$ if $W$ is
obtained by attaching a handle of index $k+1$ to $M\x [0,1].$ It
means that
 $W$ is diffeomorphic to
$M\x [0,1]\cup_f(D^{k+1}\x D^{2n-k+1}),$ where $f:S^k\x
D^{2n-k+1}\rar M\x \{ 1\}$ is the attaching map of the handle.  In
particular, $\partial W=M\cup (-M'),$ where $M'=(M-f(S^k\x
D^{2n-k+1}))\cup (D^{k+1}\x S^{2n-k})$ is the result of the surgery
on $M.$

The following classical result of Eliashberg \cite{E} (cf. also
\cite{W} and Ch. 6 of \cite{G2}) is the basic tool to construct some
examples.

\begin{theorem}\label{eli} Let $\lambda$ be a contact form on a
$(2n+1)-$dimensional manifold $M$ and let $W$ be the trace of a
surgery on $M$ of index $k+1$ with $1\leq k\leq n$ and $n>1.$ If the
almost complex structure on $M\x [0,1]$ determined by $\lambda$
extends to $W,$ then there exists an exact symplectic form $\omega$
on $W$ such that $\omega$ is the symplectization of $\lambda$ near
$M\x \{ 0\}$ as well as  the symplectization of a contact form in a
collar of $M'.$ In particular, $M'$ admits a contact form.
Furthermore, if $V$ is a compact connected  almost complex
$(2n+2)$-dimensional manifold ($n>1$) and $V$  admits a Morse
function maximal on $\partial V$ such that indices of all critical
points are less or equal to $n+1,$ then $V$ admits a symplectic
structure with convex boundary (the boundary is of contact type). A
Morse function with the required properties exists if and only if
$V$ has the homotopy type of a CW-complex of dimension at most
$n+1.$
\end{theorem}

 Let us call any manifold $V$ having the above properties
of {\it Weinstein type}. Thus the double of a manifold of Weinstein
type admits a strong symplectic fold.

\begin{remark}\label{-1}  The
contact surgery in dimension 4 requires some additional assumption
on framings of the attaching spheres of 2-handles, see \cite{GR,E}
or \cite{G2}, Ch. 6.3, 6.4.
\end{remark}

We can give now examples of whole families of  contact manifolds.

\begin{proposition} The following
manifolds admit  contact structures:

\begin{enumerate}\label{exapp}

\item \label{a1} $S^{k_1}\x ... \x S^{k_r}$
if $k_1+...+k_r$ is odd;
\item \label{a2} $S^{k_1}\x ... \x S^{k_r} \times X$  if $X$ is
a closed contact manifold  and  $k_1+...+k_r$ is even;
\item \label{a3} $M \times S^{k_1}\x ... \x S^{k_r},$ if $M$ is
a closed manifold with a strong  symplectic fold and  $k_1+...+k_r$
is odd;
\item\label{bay} $M\times X,$ if $M$ is a closed orientable 4-manifold
and $X$ is contact;
\item $\Sigma \x X,$ where $X$ is contact and $\Sigma$
is a closed oriented surface
\end{enumerate}
 \end{proposition}

{\bf Proof.} Both $D^{2k}$ and $D^{2k+1}\x S^{2l+1}$ with $k\geq l$
are Weinstein manifolds, thus taking the doubles we see that
$S^{2k}$ and $S^{2k+1}\x S^{2l+1}$  admit strong symplectic folds
with any $k,l$. Therefore the first three cases follow by induction.
To get (\ref{bay}) one has to use existence of strong symplectic
folds on closed orientable 4-manifolds  \cite{B}. In the last
statement it is enough to notice that any orientable surface has a
strong symplectic fold. This statement was first proved in \cite{Bo}
for $\Sigma$ of genus $g>0$ and for $\Sigma = S^2$ in \cite{BCS}.\bk
\medskip

Any Lie group of odd dimension is obviously almost contact. However,
no general construction of contact forms on compact Lie group is
known. It can be proved that except for rank 1 there is no
G-invariant contact forms on G. The product $S^3\x S^3\x S^3$ is an
example of simply connected Lie group  which admits a contact form
but no $G-$invariant contact form. Some examples of contact forms on
quotient spaces $G/H$ which do not admit $G-$invariant contact forms
can be obtained from Theorem \ref{mainthfib}. For instance, the
following is true.

\begin{proposition} For any even $n,$ the homogenous space
$SO(n+3)/SO(n)$ is
contact, but admits no $SO(n+3)-$invariant contact form.
\end{proposition}

{\bf Proof.} The space $SO(n+2)/SO(n)$  has a $SO(n+2)-$invariant
contact form given by the circle fibration $SO(n+2)/SO(n)\rar
SO(n+2)/(SO(n)\x SO(2))$ with symplectic base. Moreover, the space
$SO(n+3)/SO(n)$ has no $SO(n+3)$-invariant contact form. Both
statements follow from Alekseevski's description of contact
homogeneous spaces \cite{A}. Consider now the bundle
$SO(n+3)/SO(n)\rar SO(n+3)/SO(n+2)$ with fiber $SO(n+2)/SO(n).$  If
$n$ is even, then on the base $SO(n+3)/SO(n+2)= S^{n+2}$ we have the
obvious strong symplectic fold. The structure group of the bundle is
$SO(n+2),$ thus the assumptions of Theorem \ref{mainthfib} are
satisfied. In fact, one can use Corollary \ref{ctctoversph} to show
that  there exists a contact form on the total space of the bundle.
\bk

Another example of this type  is the space $SO(2k+1)/SU(k)$ of
"special unitary twistors"  on $S^{2k}$ which is fibered over
$SO(2k+1) {/}SO(2k)=S^{2k}$ with fiber $SO(2k){/}SU(k).$

We will describe now examples of a modification which can be
performed on a manifold with a strong symplectic fold. Assume that
$M^{2m}$ admits a strong symplectic fold $W_-\cup W_+$ with the fold
locus $N.$ We say that a surgery on a sphere $S^{k-1}\subset M$ is
{\it symmetric of index $k,$} if it is performed using an embedding $\phi
: S^{k-1}\x D^{2m-k+1}\rar M$ such that $\phi = \phi_0\x id_{D^1},$
where $\phi_0: S^{k-1}\x D^{2m-k}\rar N$ is an embedding and $D^1$
corresponds to the transversal disk of a tubular neighborhood  of
$N.$

\begin{proposition}\label{surgery} Consider a manifold $M$ of dimension $2m>4$
with a strong symplectic fold of convex type. If  $M'$ is obtained
from $M$ by a symmetric surgery of index $k\leq m$ such that the
stable almost complex structure of $M$ extends to $M'$, then $M'$
has a strong symplectic fold structure.
\end{proposition}

{\bf Proof.} We have $M'=(M-\phi (S^{k-1}\x D^{2m-k+1}))\cup (D^k\x
S^{2m-k}).$ Decompose $S^{2m-k}$ into the sum of two disks $D_-\cup
D_+$ such that the decomposition corresponds to cutting the sphere
by $N.$ Then we obtain a decomposition $M'=W'_-\cup W'_+$ such that
$W'_{\pm}= W_{\pm} \cup (D^k\x D_{\pm}).$ Thus both parts are given
by attaching handles  $D^k\x D_{\pm}$ of index $k$ to respectively
$W_-,W_+.$ Because $k\leq m$ and by assumption the almost complex structures
on $W_-, W_+$ extend to these handles, given symplectic  forms
extend to $W_{\pm}.$ \bk

\begin{corollary} If $M^{2m}$ admits a strong symplectic fold ($k+n=2m$), then
so does the connected sum $M\# (S^k\x S^n). $
\end{corollary}

{\bf Proof.} The  proposition can be applied, since connected sum
with $S^k\x S^n$  is obtained by the surgery on a trivially embedded
sphere $S^{k-1}$ (or on $S^{n-1}$) and thus we can assume that
$k\leq n.$ \bk

\section{Concave folds and strong symplectic folds of general
type}\label{conc}

Till now we considered decompositions of a manifold $M$ into the sum
of two exact symplectic cobordisms $W_1$ and $W_2$ having the same
contact boundary $N$ at their convex ends. If $M$ is closed, the
symplectic cobordisms cannot have concave ends, thus they should be
symplectic fillings of the contact form on the fold locus. Our
present purpose is to extend this construction  to the case where
$M$ is decomposed into several pieces, each being a symplectic
cobordism having possibly concave ends as well.

In case when $W_1,W_2$ meet in such a way that one of the ends is
concave and one is convex (and the contact forms at the boundary are
equal), one can apply standard gluing of two symplectic cobordisms,
which assembles two symplectic cobordisms into one, simplifying the
decomposition. Thus the substantial cases are when two convex ends
or two concave ends meet.

Consider now the case of concave ends of two symplectic cobordisms
meeting  at $(N,\la ).$ We assume that in a collar neighborhood $N\x
[-1,1]$ of $N$
 we have the form $e^{-t}\la$ for $t \in [-1;0]$ and
$e^{t}\hat\la$ for $t \in [0;1]$. Note that the orientations given
by  the forms on the two sides coincide, unlike the case of convex
folds.

We explain now how to use Theorem \ref{mainth} to obtain a contact
form on the product of the sum of such two cobordisms by a contact
manifold $X.$

\begin{lemma}\label{swap} Suppose that $(X,\al)$ and $(N,\la )$ are
closed contact manifolds. Then there exists a contact form on $X\x N
\x [-1,1]$ equal to $\al + e^{-t}\la $ near $X \x N \x \{-1\}$ and
to $\al + e^t\hat\la$  near $X \x N \x \{ 1\} .$

\end{lemma}

{\bf Proof.} We apply Theorem \ref{mainth} after switching the role
of $X$ and $N.$

For this purpose we define positive functions

$$\label{g1}g_1(t)=\left\{ \begin{array}{rl}
1 &\textrm{near $t=-1$}
\\ e^t & \textrm{near $t=-\frac 12$}
\end{array} \right.,$$

\NI on $[-1,-\frac12]$ and

$$\label{g2}g_2(t)=\left\{ \begin{array}{rl}
e^{-t} & \textrm{near $t=\frac 12$}
\\ 1 & \textrm{near $t=1$}
\end{array} \right..$$

\NI on $[\frac12,1].$

\NI The contact form $g_1(t)(\al + e^{-t}\la)$ on $X\x N\x
[-1,-\frac12]$ extends  to $e^t\al + \la$ on $X\x N\x [-\frac12,0].$
Similarly, the contact form $g_2(t)(\al + e^{t}\hat\la)$ on $X\x N\x
[\frac12,1]$ extends  to $e^{-t}\al + \hat\la$ on $X\x N\x [0,
\frac12].$ Thus we can apply Theorem \ref{mainth} to construct the
form with required properties.

\bk

In Lemma \ref{swap}, in order to get a contact structure on the
product of this manifold by a contact one, we need a contact form
$\la $ on one end of $N\x [-1,1]$ and $\hat\la $ on the other, while
on the contact factor the form does not depend on $t\in [-1,1].$
This is too restrictive for applications and we will show that the
construction of contact forms is possible also if the  pair of forms
is  $(\la + \al)$ and $(\la + \hat\al )$ on the two sides.

\begin{lemma} \label{concaveends}   Suppose that
$(X,\al)$ and $(N,\la )$ are closed contact manifolds. Then there
exist two contact forms on $X\x N \x [-1,1],$ both equal to
$e^{-t}\la + \al$ near $N \x \{-1\}$ and one equal to $e^t\hat\la +
\al ,$ the other to $e^t\la + \hat\al$ near $N \x \{ 1\} .$
\end{lemma}

{\bf Proof.} In case of the pair $(e^{-t}\la + \al , e^{-t}\hat\la +
\al )$ we apply Lemma \ref{swap}. If we have the pair $(e^{-t}\la +
\al , e^{-t}\la + \hat{\al})$ we use Theorem \ref{mainth} and Lemma
\ref{swap} to $X\x M\x [-1,1]$ divided into 4 parts:

\begin{enumerate}

\item $X\x N\x [-1,-\frac12]$ with the form $e^{-t}\la +\al ,$

\item $X\x N\x [-\frac12,0]$ with the form $e^{t}\hat\la +\al ,$

\item $X\x N\x [0,\frac12]$ with the form $e^{-t}\hat\la +\hat\al ,$

\item $X\x N\x [\frac12,1]$ with the form $e^{t}\la +\hat\al .$

\end{enumerate}

The forms are defined such that crossing the convex fold at $0$
corresponds to passing from $(\al ,\la )$ to $(\hat\al ,\la )$ and
crossing concave folds at $\pm \frac12$ is the swap between $(\al
,\la ), (\al ,\hat\la )$ and back. In all cases one of the
previously described constructions works. Thus we get a contact form
on  $X\x N.$

\smallskip

\centerline{\includegraphics[width=8cm,height=4cm]{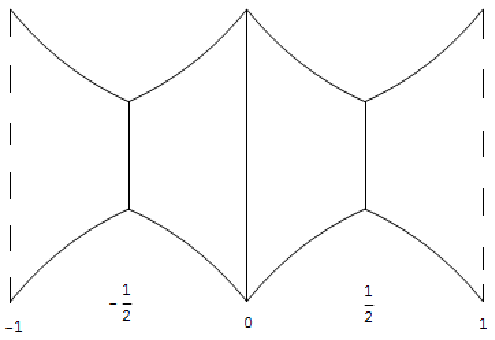}}

\bk

Now we are in position  to extend the notion of strong symplectic
fold to allow concave folds.

Consider a closed hypersurface $N\subset Int\, M$ and denote by
$W_i, i=1,..,k$ the connected components of $M-N$ compactified by
adding adjacent components of $N.$ Hence $W_i$ is just the  closure
of a component of $M-N.$  Let $N=\bigcup_sN_s$ denote the
decomposition of $N$ into the sum of connected components.

\begin{definition} A strong symplectic fold on a compact
 manifold $M$ is given by:

\begin{enumerate}

\item  a decomposition $\{W_i\}_{i \in I}$ of compact codimension 0 submanifolds,
$M = \bigcup_iW_i,$ obtained by cutting $M$ by a hypersurface
$N\subset Int\, M;$

\item a family of contact forms
$(N_s, \la_s);$

\item exact symplectic forms $\omega_i=d\beta_i$ on $W_i$
such that each $\omega_i$ yields a symplectic cobordism structure on
$W_i$ with some convex ends and some concave ends and each pair
$\omega_i, \omega_j$ satisfies one of  the following compatibility
condition for every connected component $N_s$  of $N$ with
$N_s\subset W_i\cap W_j$:

\begin{enumerate}

\item\label{ssf21}  $\beta_i = e^{t}\la_s $ in
$N_s\times [-1,0]\subset W_i$ and $\beta_j=e^{-t}\la_s$ in
$N_s\times [0,1]\subset W_j$  where $t$ is the parameter of $[-1,1]$
(convex fold: a convex end of $W_i$ meets a convex end of $W_j$ at
$N_s$);

\item\label{ssf22} $\beta_i = e^{-t}\la_s$ on $N_s\times
[-1,0]\subset W_i$ and $\beta_j= e^{t}\la_s$ on $N_s\times
[0,1]\subset W_j,$ where $t$ is the  parameter of $[-1,1]$ (concave
fold: a concave end of $W_i$ meets a concave end of $W_j$);

\end{enumerate}

\end{enumerate}

\end{definition}

\NI As before, the hypersurface $N$ is called the {\it fold locus.}
We assumed that every piece $W_i$ is a symplectic cobordism, hence
the forms $\omega_i$ are either convex or concave along any
component of the boundary of $M$.

From the discussion of this section we obtain the following
extension of Theorem \ref{mainth}.

\begin{theorem}\label{mainthgen} If $X$ is a
closed contact manifold and $M$ admits a strong symplectic fold,
then $X\x M$ is contact.
\end{theorem}   \bk

\section{Some further applications}\label{furapp}

To illustrate usefulness of concave folds consider the question of
fillability of contact manifolds (by a symplectic one). It is
well-known that no overtwisted contact form $\la$ on a compact
3-manifold $M$ is fillable, i.e., there is no compact manifold with
boundary of contact type (convex boundary) having overtwisted
contact form on the boundary. Constructions based on fibrations, for
instance the open book technique, lead to the following question. Is
there a similar obstruction to fill up by a compact contact manifold
the product of an overtwisted 3-manifold by $S^1?$ In other words,
we ask if the form $e^{-t}\la + d\theta$  on $M\x [0,\eps )\x S^1\ $
can be extended to a contact form on a compact manifold $W$ such
that $\partial W = M\x \{ 0\}\x S^1$ ($d\theta$ denotes the
standard form on $S^1$). Below we show examples that
fillability in this sense is possible.

Given two connected contact manifolds $(X,\al ), (X', \al')$
oriented compatibly with contact structures,  one can perform
1-surgery such that the resulting manifold is the connected sum $X\#
X'.$ Then by the contact surgery (Theorem \ref{eli}) we get a
contact form on the connected sum. Since we need some choices to
perform such operation, the result is not defined uniquely, but its
homotopy class is already unique. By slight abuse of language we
denote the contact form obtained in this way by $\al\#\al' .$

\begin{proposition}\label{exttodisk}
If $n>0$ and $\lambda $ is any contact form on $S^{2n+1},$ then the
form $e^{-t}(\lambda \# \hat\lambda ) + d\theta$ on a collar of the
boundary  $(S^{2n+1}\#S^{2n+1})\x [0,\eps )\x S^1 \subset D^{2n+2}\x
S^1$  extends to a contact form on $D^{2n+2}\x S^1.$
\end{proposition}

{\bf Proof.} Consider the symplectizations $e^{-t}\lambda$ on
$S^{2n+1}\x [-1,0]$ and $e^t\hat\lambda$ on $S^{2n+1}\x [0,1].$
Gluing these manifolds along $S^{2n+1}\x \{ 0\}$ we get a manifold
with a concave fold $S^{2n+1}\x \{ 0\}$ and boundary $S^{2n+1}\cup
-S^{2n+1}.$ We can perform contact 1-surgery by adding a 1-handle to
the boundary which makes the boundary connected and diffeomorphic to
$S^{2n+1}.$ The manifold $W$ obtained by the surgery is
diffeomorphic to $S^{2n+1}\x S^1 - D^{2n+2}.$ Using Theorem
\ref{eli} for this handle we get a strong symplectic fold on $W$
extending the symplectizations and with the boundary $(S^{2n+1},
\lambda \# \hat\lambda )$ of contact type (note that we still have
the fold $S^{2n+1}\x \{ 0\}$ in the interior of $W).$ By Lemma
\ref{swap}, there is a contact form $\eta$ on $W\x S^1$ equal to
$\lambda\#\hat\lambda + d\theta$ on the boundary. Denote by
$S\subset Int\, W$ the circle given as the sum of intervals $x_0\x
[-1,1]\subset S^{2n+1}\x [-1,1]$ and $y_0\x [-1,1]$ in the handle,
where $y_0\x \{ \pm 1\}$ are attached to $x_0\x \{\pm1\}$ by the
attaching map of the handle. The (topological) surgery of index 2 on
$W$ with the attaching circle $S$ and the standard framing of the
normal bundle yields the disk $D^{2n+2}.$ Moreover, the standard
almost complex structure on $W$ extends to the 2-handle. To finish
the proof we have to show that the surgery applied to $\eta$ yields
another contact form on its result.  $D^{2n+2}\x S^1$ is obtained
from $W\x S^1$ by the 2-surgery multiplied by $S^1.$ The product of
a 2-handle  by $S^1$ decomposes into two handles on $W\x S^1$, one
of index 2 on $W\x S^1$  and one of index 3 attached to the result
of the first surgery. Since the manifold $W\x S^1$ is of dimension
at least $5$ and the given almost contact structure is compatible
with the surgeries, we get a contact form on $D^{2n+2}\x S^1.$
Finally, the surgeries are done in the interiors of manifolds in
each step of the construction, hence they preserve the form we have
obtained previously in a neighborhood of the boundary sphere.

 \bk

Thus we get the following corollary that was first proved in \cite{EP}.

\begin{corollary}\label{extover}  There exists an overtwisted contact
form $\lambda$ on $S^3$ such that the form  $e^{-t}\lambda +
d\theta$ extends from a collar $S^3\x S^1\x [0,\eps )$ to a contact
form on $D^4\x S^1.$
\end{corollary}

This property can be applied to prove the following  special case of
results proved in \cite{EJ,CPP}.

\begin{proposition}\label{4dim} If $M^5$ is closed almost
contact and admits  an
open book decomposition with trivial monodromy, then it is contact.
\end{proposition}

{\bf Proof.}  Let $P$ denote the page of the open book. The almost
contact structure of $M$ gives a stably almost complex structure on
$P.$ For an open manifold stably almost complex structure determines
an almost complex structure. It follows from basic facts of the
Morse - Smale theory that there exists a Morse function $f:P\rar
[0,4]$ with one minimum $(=0)$, constant and  maximal $(=4)$ on
$\partial P.$ This function has critical points only of indices
$q=0,1,2,3$ and such that the value of $f$ at a critical point of
index $q$ is $q.$ Denote $W_i= f^{-1}[i-\frac12,i+\frac12],\ i=
0,1,2,3,4.$   Then $W_0$ is diffeomorphic to $D^4,$ $W_i$ contains
only critical points of indices $i$ and $W_4=\partial P\x
[\frac72,4].$  Let $\lambda$ be an overtwisted contact form on $S^3$
such that $e^{t}\lambda + d\theta$ extends to a contact form on
$D^4\x S^1.$ Since $P$ is almost complex, then by the contact
surgery Theorem \ref{eli} we extend the form $e^{t}\lambda$ to 1-handles of
$W_1.$ This makes $W_1$ a symplectic cobordism with concave end
$f^{-1}(\frac12)$ and convex end $f^{-1}(\frac32).$ Since the
surgeries can be performed far from overtwisted disks, the contact
form on the latter can be assumed again overtwisted. On an
overtwisted 3-manifold one can perform contact surgery on every
framing, so this holds for $W_2.$  In the same manner we make $W_3$
a symplectic cobordism with concave end $f^{-1}(\frac72)$ and convex
end $f^{-1}(\frac52).$ Namely, we use the (unique up to homotopy)
overtwisted form $\mu$ representing the   almost contact structure
of $f^{-1}(\frac72).$ In this way we get symplectic structures which
agree with the almost complex structure of $P.$ Since the homotopy
class of an overtwisted form is determined by the homotopy class of
the contact distribution, the contact forms on $f^{-1}(\frac52)$
obtained from $W_2$ and $W_3$ are homotopic, hence by Remark
\ref{concordance} can be assumed equal. Finally, on $W_4$ we put
symplectization of the form $\hat\mu , $ where $\mu$ is the form
used in $W_3.$ In this way we get a strong symplectic fold on
$W_1\cup W_2\cup W_3\cup W_4$ with fold locus $f^{-1}(\frac52)\cup
f^{-1}(\frac72),$ where the fold at $f^{-1}(\frac52)$ is convex  and
at $f^{-1}(\frac12)$ is concave. Therefore, by Theorem \ref{mainth}
and Lemma \ref{swap} we have a contact form on the product with
$S^1.$ Since the form on $f^{-1}(\frac12)\x S^1$ extends to $D^4\x
S^1,$ we get also a contact form on $P\x S^1.$ By the construction,
in a collar of $\p P\x S^1$ this form is the product of a convex
form on $W_4$ by the standard form on $S^1$ at $\partial P.$ It can
be extended to $\p P\x D^2\subset M$ exactly as it is done in the
case of the open book construction. This completes the proof. \bk

\ms

Let us illustrate Theorem \ref{mainthgen} by the following examples.

\begin{example}  If $M$ is a $S^1$-bundle over $X\x N$ with
contact $(X,\la),(N,\la'),$ then $M$ is contact (in particular,
$X\x S^1\x S^1$ is). Begin with the trivial bundle. Write $N\x S^1 =
N\x [0,\frac14]\cup N\x [\frac14,\frac12]\cup N\x
[\frac12,\frac34]\cup N\x [\frac34,1],$ where $N\x \{ 0\}$  and $N\x
\{ 1\}$ are identified. On these four parts put $e^t\lambda',
e^{-t+\frac12}\lambda', e^{t-\frac12}\hat\lambda',
e^{-t+1}\hat\lambda',$ respectively. This gives a strong symplectic
fold structure on $N\x S^1.$  Now take products with $(X,\la)$ for
$N\x [0,\frac14], N\x [\frac34,1]$ and with $(X,\hat{\la})$ for $N\x
[\frac14,\frac12], N\x [\frac12,\frac34].$ So we have the following
sequence of forms:

$$e^t\la' +\la ,\ \  e^{-t+\frac12}\la'+\hat{\la},\ \
e^{t-\frac12}\hat{\la'}+\hat{\la},\ \   e^{-t+1}\hat{\la'}+\la
.$$

\centerline{\includegraphics[width=8cm,height=4cm]{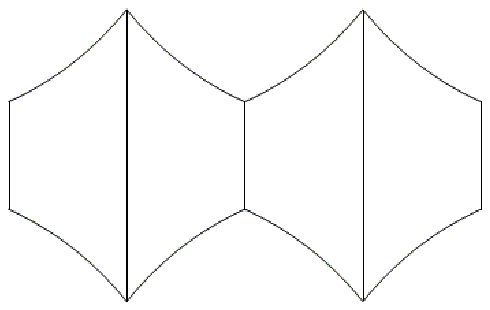}}

By Theorem \ref{mainthgen}, there exists a contact form  on $X\x N\x
S^1.$ By \cite{G3}, this extends to any circle bundle over $X\x N.$

\end{example}

\begin{example}\label{blow} Consider a closed contact manifold
$M$ of dimension $5$ and a homotopically trivial circle $S$ embedded
in $M.$ Then the manifold $M'$ obtained from $M$ by the blow-up
along $S$ is contact.
 \end{example}

{\bf Proof.} We can deform the given contact form on $M$ to one
given by an open book with the binding $B,$ the page $P$ and the
fibration $E\rar S^1.$ Then $S$ can be deformed to a section of the
fibration, say to a circle given by a point near $\p P,$ where the
fibration is product. A tubular neighborhood of such $S$ is the
product of a small disk $D^4_0\subset Int\, P$ by $S.$ On
$\overline{\mathbb CP^2}$ there is a strong symplectic fold $W_-\cup
W_+$ of convex type by \cite{B}. Cutting another  small (Darboux)
disks $D^4_1$ in $W_-$ and identifying boundary spheres of $D^4_0,
D^4_1$  we get the connected sum $P \# \overline{\mathbb CP}^2$ and
a strong symplectic fold on it (with concave fold at the connected
sum sphere). Consider the following decomposition of $M:$ $B\x D^2,$
the product neighborhood of the binding, $(U-D^4_0)\x S^1,$ where
$U$ is a collar of $\p P, \ (W_--D^4)\x S^1, W_+\x S^1$ and the
fibration over $S^1$ with fiber $P-U$ given by the open book
structure. On the fibration we have a contact form. By the assumption, this form is
product near the boundary. Other pieces are products, thus we can
apply Theorem \ref {mainthgen} to get a contact form on $M.$

\bk

\begin{remark}  Note that $P\# \overline{\mathbb CP}$
does not admit any exact symplectic form with contact type boundary,
so the example cannot be obtained by modification of the open book.
It was explained to us by Andr\'as Stipsicz that this  property
follows from the fact that any spherical homology 2-class in a
closed 4-manifold with self-intersection number $-1$ is represented
by a symplectic submanifold. The same argument, combined with a
result of McDuff \cite{McD} shows that there is no strong symplectic
fold of convex type on $\mathbb CP^2-Int\, D^4.$
 \end{remark}

\section{A generalization of the open book construction}\label{genopen}

We give now a generalization of the open book construction allowing
bindings of codimensions greater that 2. This is a decomposition of
a manifold into two pieces. This decomposition is much more symmetric than the
open book.

   Consider two compact manifolds $X,Y$  with non-empty
boundaries, of dimensions $2n,2m$ respectively. Assume that they are
endowed with  exact symplectic forms $\omega_X=d\beta_X,\
\omega_Y=d\beta_Y,$ both with convex type boundaries.  Let
$\beta_X=e^s\mu_{\p X},\  \beta_Y=e^{-s}\mu_{\p Y}$ in collars $\p
X\x (-1,0], \p Y\x [0,1)$ of boundaries, where $\mu_{\p X}, \mu_{\p
Y}$ are some contact forms. In this notation $s\in [-1,1]$ and both
boundaries correspond to $s=0.$ Let $E$ be the total space of a
bundle over $\partial Y$ with fiber $X$ defined on a hypersurface
$H_{\p Y}\subset
\partial Y$ and with the structure group $Ex(X,\partial X,
\beta_X)$  of exact symplectomorphisms of $X$ equal to the identity
near the boundary.  Similarly, assume that the bundle $F\rar\partial
X$ with fiber $Y$ is defined on a hypersurface $H_{\p X}\subset
\partial X,$ and its structure group is $Ex(Y,\partial Y, \beta_Y).$
The assumptions on structure groups imply that $\partial E =
\partial X\x \partial Y = \partial F.$

\begin{proposition}
Under the above assumptions, $E\cup_{\p X\x \p Y} F$ is contact.

\end{proposition}

{\bf Proof.} Consider  $\tilde X = X\cup \p X\x [0,\log R_X]$
obtained from $X$  by adding a long collar,  with $\beta_X =
e^s\mu_{\p X}$ for $s\in [-1,\log R_X].$ In this way the contact
form on the boundary is multiplied by the constant $R_X.$
Analogously, $Y$ is enlarged to $\tilde Y = Y\cup \p Y\x [-\log
R_Y,0]$ with $\beta_Y = e^{-s}\mu_{\p Y}$ for $s\in [-\log R_Y,1]$
(we assume $R_X,R_Y\geq 1).$ Let $\tilde E$ denote the obvious
extension of $E$ to a bundle with fiber $\tilde X,$ and similarly
$\tilde F$ the extension of $F.$ Proposition \ref{exactbundles}
gives  a contact form on $\tilde E$ equal to $R_Y\mu_{\p Y}+\beta_X$
near $\p \tilde E = \p Y\x \p X.$ The choice of $R_Y$ which yields
contactness is determined by the behavior of the forms in the
tubular neighborhood of $H_{\p Y}.$ We claim that the choice depends
only on $X$ (not on $\tilde X$ and $R_X).$ To see this, let us
calculate $\eta_E d\eta_E^{n+m-1}$ for $\eta_E = R_Y\mu_{\p Y} +
\beta_X + ud\psi$ in $\p X\x [-1, \log R_X]\x H\x [-1,1].$ Since
$\beta_X=e^s\mu_{\partial X}$ for $ \ s\in [-1,\log R_X],$

$$\eta_E d\eta_E^{n+m-1}= e^{ns}\mu_{\partial X}
d\mu_{\partial X}^{n-1} \left( D_1R_Y^m\mu_{\partial Y}
d\mu_{\partial Y}^{m-1}ds+ D_2 R^{m-1} u'\mu_{\partial Y}
d\mu_{\partial Y}^{m-2} dsdt d\widetilde{\psi}\,+ \right.$$
$$\left. D_3 R^{m-2} u' d\mu_{\partial Y}^{m-1}dt d\widetilde{\psi}+
D_4R^{m-1} d\mu_{\partial Y}^{m-1} ds d\widetilde{\psi} \right)
,$$ where $D_i,i=1,2,3,4$ are again constants depending only on
$m,n.$

It follows from this formula that the choice of $R_Y$ is independent
of the extension by the long collar and our claim follows. Thus we
can choose $R=R_X=R_Y$ such that there are contact forms on $\tilde
E$ and $\tilde F$ that restrict to $R(\mu_{\p X} + \mu_{\p Y})$ on
$\p E = \p F = \p X\x \p Y.$ Let $K=\log R.$ After the change the
parameter in $[-K,1]$ replacing $s$ with $s+2K,$ the form $\eta_F$
on $\p X\x \p Y\x [K,1+2K]$ becomes  $R\mu_{\p X} + e^{-s+K}\mu_{\p
Y}.$

Let $\psi:[K,1+2K]\rar \mathbb R$ be a positive smooth function such
that $\psi=e^{s-K}$ near $s=K$ and $\psi=1$ near $1+2K,$ regarded as
a function on $F$ (we simply extend it from the collar $\p F\x
[0,1]$ to whole $F$). Then  $\psi\eta_F$ is contact and it smoothly
agrees with $\eta_E$ along $\partial E = \p F = \p X\x \p Y.$  Thus
we get a smooth contact form on $E\cup F.$ \bk

\section{Concluding remarks}\label{finalsec}

We do not know any example of closed stably almost complex manifold
which admits no strong symplectic folds.  On the other hand, it is
anything but obvious if any symplectic manifold has a strong
symplectic fold. In particular, it would be interesting to decide
whether complex projective spaces admit strong symplectic folds.

The standard Morse - Smale theory shows that for any closed manifold
$M^{2m}$ one can find a decomposition $M=W_+\cup_NW_-,$ where $N =
\partial W_+=\partial W_- = W_+\cap W_-$ with both $W_+,W_-$ having
the homotopy type of complexes of dimension at most $m.$ If $M$ is
stably almost complex, then  $W_{\pm}$ are almost complex and we
have exact symplectic forms on both parts by contact surgery. The
resulting contact forms $\lambda_-, \lambda_+$ on $N$ define
homotopic almost contact structures on $N,$ but the question whether
they are homotopic (as contact forms) is apparently difficult. If
they do, we would get a strong symplectic fold of convex type on
$M.$ Quite possibly, the general type of strong symplectic folds can
be useful in this problem, as the arguments used  for Example \ref{blow} and
Proposition \ref{4dim} indicate.

\ms

\section{Appendix A: Contact piecewise fibered
structures}\label{locfib}

In this appendix we describe a preliminary version of a structure
generalizing all cases we considered till now and still sufficient
to provide a contact form on a manifold endowed with such structure.
This is obtained by localization, requiring that each piece of such
decomposition is one of described previously with appropriate
compatibility conditions along  intersections assumed.

Let $Y$ be a compact orientable manifold. Given  a hypersurface
$H\subset Int Y,$ let $\{ Y_{i}\}$ denote the collection  of
connected components of $Y-H$ compactified by adding components of
$H$ contained in the closure of $Y_i.$ Our basic assumption is that
each $Y_{i}$ is a fibration of one of the following two types:

 \begin{enumerate}

 \item  a contact fibration with a closed contact fiber $(X_i,\al_i)$
over an exact symplectic cobordism $(W_i, d\mu_i)$ trivial in a
neighborhood of $\p W_i, $ or

\item  the fibration over a closed contact manifold $(X_i,\al_i),$
defined on a hypersurface in $X_i,$ such that the fiber is an exact
symplectic cobordism $(W_i,d\mu_i)$ and the structure group is  the
group $Ex(W_i,\p W_i,\mu_i)$ of exact symplectomorphisms equal to
the identity in a collar of $\p W_i.$

\end{enumerate}

If this is satisfied, then every component of $H$ is the product of
$X_i$ by  a component of the boundary of the symplectic cobordism
$W_i.$  Let us denote by $N_{is}, s=1,..,l_s$ components of $\p W_i$
and by $ \la_{is}$ the contact form induced on $N_{is}$ by $\mu_i$
(which is either convex or concave at $N_{is}).$

If $N_{is}=N_{jr}$ is a  connected component of the intersection
$Y_i\cap Y_j\cap H,$ then we assume that one of the following
conditions is satisfied:

\begin{enumerate}
\item \label{convfold} $N_{is}$ is a convex end of $W_j,$
$N_{jr}$ is a convex end of $W_j$ and $X_i = X_{j};$

\item  \label{concfold} $N_{is}$ is a concave end of $W_i,$
$N_{jr}$ is a concave end of $W_j$ and $X_i=X_j;$

\item \label{girfold} $N_{is}$ is a convex end of $W_i,$ $N_{jr}$
is a convex end of $W_j,$  $N_{is}=X_{j},\  N_{jr}= X_i;$

 \item \label{girfold1}  $N_{is}$ is a concave end of $W_i,$ $N_{jr}$
is a concave end of $W_j,$  $N_{is}=X_{j},\  N_{jr}= X_i;$

 \item \label{mix1} $N_{is}$ is a concave end of $W_i,$
$N_{jr}$ is a convex end of $W_j$ and $X_i = X_j;$

\item   \label{mix2} $N_{is}$ is a convex end of $W_i,$
$N_{jr}$ is a concave end of $W_j$ and $X_i = X_j.$

\end{enumerate}

Finally, we assume compatibility  of the forms on the adjacent ends
of $Y_i's.$ In all the cases above we require one the following
conditions, according to the list above:

\begin{enumerate}

\item $\la_{is} = \la_{jr}$ and $\al_{i}=\hat\al_{j}$
or $\la_{is} = \hat\la_{jr}$ and $\al_i=\al_j;$

\item  $\la_{is} = \la_{jr}$ and $\al_i=\hat\al_j$
or $\la_{is} = \hat\la_{jr}$ and $\al_i=\al_j;$

\item $\la_{is} = \al_{j}$ and $\al_{i}=\la_{jr};$

\item $\la_{is} = \al_j$ and $\al_i=\la_{jr};$

\item $\la_{is} = \la_{jr}$ and $\al_i=\al_j;$

\item $\la_{is} = \la_{jr}$ and $\al_i=\al_j.$

\end{enumerate}

If $\p Y_i$ contains a connected component of $\p Y,$ then in a
collar of that component we have the product of $X_i$ and an end of
$W_i$ (either convex or concave).

\begin{remark}  We allow a component of $H$ to be the boundary
of two different ends of one $Y_i$ (when $i=j$ in the list above).
In particular, it is possible that $Y-H$ is connected.\end{remark}

 One can explain our assumptions by saying that the fold locus $H$
divides the manifold $M$ into a number of fibrations carrying
contact fibered structure with both fibrations and forms product
near any component of $H.$  Under our compatibility conditions we
can apply either Theorem \ref{mainthfib} or Lemma \ref{concaveends}.

\begin{definition} A decomposition of $M$ satisfying the
assumptions above is called {\it a contact piecewise  fibered
structure} on $M.$
\end{definition}

\begin{theorem}\label{comp} If $M$ admits a contact piecewise  fibered
structure, then $M$ is contact.
\end{theorem}

{\bf Sketch of the proof.}  Consider a component $Y_{i}$ of the
decomposition. As we explained in Sections \ref{prel} and
\ref{mainsec}, it admits a contact form equal to $\la^{\eps}_{ij} +
p^*\al_i,$ or to $p^*\la_{ij} + \al_i,$ in a collar of the j-th
component of $\p Y_i,$ depending on the type of the fibration on
$Y_i.$ Furthermore, $\eps=\pm 1$ depending on convex/concave type of
the fold. Under the compatibility conditions we use Theorem
\ref{mainthfib} or Lemma \ref{concaveends} to extend those forms
through $H$ and we get a global contact form on $M.$

\bk

\begin{remark}  One can allow that instead of equalities in the
compatibility conditions one assumes equality up to homotopy, for
instance up to the multiplication by a constant. This can be always
reduced to the equality case  by extending the adjacent end (which
is $e^{\pm t}\la_j,\ t\in [0,1])$ from $[0,1]$ to $[0,R]$ for $R$
appropriately chosen and  applying the trick of  Lemma
\ref{concaveends}.
\end{remark}

\section{Appendix B: Computations in Mathematica}

We present here some of the calculations which led us to the proof
of Theorem \ref{mainth}. The result was first checked using
Mathematica's package ''Differential forms'' (Frank Zizza and Ulrich
Jentschura \cite{FZ}) in low dimensions. Namely, for $t1=0$ (for
technical reasons we slightly change notation to adapt it for our
purposes) and around a point $(b,d,n,0) \in B \times D^2 \times N
\times I$ we take coordinate system in which $\beta=d[z1]+x1 d[y1],
\lambda = d[z2]+x2 d[y2].$ Further, on disk $D^2$ we take coordinate
system $(x,y).$ In these coordinates we set $h_1=2-(x^2+y^2)^2$ and
$h_2=x^2+y^2$ (hence in the formula below $h_1$ is equal to $2-r^4$
near $r=0$ so that it is of class $C^3$). Then the following
expressions are equal respectively to $\eta$ and $d\eta:$

\medskip

\noindent\({\text{eta1}\text{:=}(2-(x{}^{\wedge}2+y{}^{\wedge}2)
{}^{\wedge}2)(d[\text{z1}]+\text{x1} d[\text{y1}])+(d[\text{z2}]+\text{x2}
d[\text{y2}])
}\)

\medskip

\noindent\({\text{deta1}\text{:=}\left(-4 x^3-4 x y^2\right) d[x]
\wedge d[\text{z1}]+\left(-4 x^3 \text{x1}-4 x \text{x1} y^2\right) d[x]\wedge
d[\text{y1}]+}\\
{\left(-4 x^2 y-4 y^3\right)d[y]\wedge d[\text{z1}]+\left(-4 x^2
\text{x1} y-4 \text{x1} y^3\right)\text{x1} d[y]\wedge d[\text{y1}]+}\\
{\left(2-x^4-2 x^2 y^2-y^4\right)d[\text{t1}]\wedge d[\text{z1}]+
\left(2-x^4-2 x^2 y^2-y^4\right)d[\text{x1}]\wedge d[\text{y1}]+}\\
{\left(2 \text{x1}-x^4 \text{x1}-2 x^2 \text{x1} y^2-\text{x1}
y^4\right)d[\text{t1}]\wedge d[\text{y1}]+d[\text{x2}]\wedge d[\text{y2}]+}\\
{d[\text{t1}]\wedge d[\text{z2}]+\text{x2}\text{  }d[\text{t1}]
\wedge d[\text{y2}]-x d[\text{t1}]\wedge d[y]+y d[\text{t1}]\wedge d[x]}\)

\medskip

\NI Now $\tau=\star(\eta \wedge (d\eta)^3)$ can be computed in two steps:
first we calculate

\medskip

\noindent\({\text{ExteriorProduct}[\text{eta1},\text{deta1},\text{deta1},\text{deta1}]}\)

\medskip

\NI and later

\medskip

\noindent\({\text{HodgeStar}[\text{\%},t[\text{x1},\text{x1}]+t[\text{y1},\text{y1}]+
t[\text{z1},\text{z1}]+t[\text{x2},\text{x2}]+}\\
{t[\text{y2},\text{y2}]+t[\text{z2},\text{z2}]+t[x,x]+t[y,y]+t[\text{t1},\text{t1}]]}\)

\medskip

\NI where the percent sign refers to $\eta \wedge (d\eta)^3.$

\medskip

\NI Then $\tau=\star(\eta \wedge (d\eta)^3)$ is given by

\medskip

\noindent $(24 \left(x^2+y^2\right)^2\text{) }d\text{x1}\text{ ${}^{\wedge}$
}d\text{y1}+(-96 x (-1+\text{x1}) \text{x1} y \left(x^2+y^2\right)^2\text{)
}d\text{t1}\text{ ${}^{\wedge}$ }d\text{x1}+$

\noindent $(-24 \text{x1} \left(x^2+y^2\right) \left(x^2+\text{x1}
y^2\right)\text{)
}d\text{x1}\text{ ${}^{\wedge}$
}d\text{z1}+(6 x \left(-2+x^4+2 x^2 y^2+y^4\right)\text{) }dx\text{
${}^{\wedge}$ }d\text{z1}$

\noindent $+(6 y \left(-2+x^4+2 x^2 y^2+y^4\right)\text{) }dy\text{
${}^{\wedge}$ }d\text{z1}+(-24 \left(x^2+y^2\right)^2 \left(-2+x^4+2 x^2
y^2+y^4\right)\text{) }
d\text{x2}\text{ ${}^{\wedge}$ }d\text{y2}+$

\noindent $(24 \text{x2}
\left(x^2+y^2\right)^2 \left(-2+x^4+2 x^2 y^2+y^4\right)\text{) }
d\text{x2}\text{ ${}^{\wedge}$ }d\text{z2}+(6 x \left(-2+x^4+2 x^2
y^2+y^4\right)^2\text{)
}dx\text{ ${}^{\wedge}$ }d\text{z2}+$

\noindent $(6 y \left(-2+x^4+2 x^2 y^2+y^4\right)^2\text{) }dy\text{
${}^{\wedge}$ }d\text{z2}+
(24 x \left(x^2+y^2\right) \left(-2+x^4+2 x^2 y^2+y^4\right)\text{)
}d\text{t1}\text{ ${}^{\wedge}$ }dy+$

\noindent $(-24 y \left(x^2+y^2\right) \left(-2+x^4+2 x^2
y^2+y^4\right)\text{) }d\text{t1}\text{ ${}^{\wedge}$ }dx$

\noindent $+(-24 (-1+\text{x1}) \text{x1} y^2 \left(x^2+y^2\right)
\left(-2+x^4+2 x^2
y^2+y^4\right)\text{) }d\text{x1}\text{ ${}^{\wedge}$ }d\text{z2}$

\medskip

\NI and $\tau^4$ is equal to
\medskip

\noindent\((-1990656 \left(x^2+y^2\right)^6 \left(-2+x^4+2 x^2
y^2+y^4\right)^3\text{)
}d\text{t1}\text{ ${}^{\wedge}$ }dx\text{ ${}^{\wedge}$ }d\text{x1}\text{
${}^{\wedge}$ }d\text{x2}\text{ ${}^{\wedge}$ }dy\text{ ${}^{\wedge}$
}d\text{y1}\text{ ${}^{\wedge}$ }d\text{y2}\text{ ${}^{\wedge}$
}d\text{z1}+(-1990656
\text{x2} \left(x^2+y^2\right)^6 \left(-2+x^4+2 x^2 y^2+y^4\right)^3\text{)
}d\text{t1}\text{ ${}^{\wedge}$ }dx\text{ ${}^{\wedge}$ }d\text{x1}\text{
${}^{\wedge}$ }d\text{x2}\text{ ${}^{\wedge}$ }dy\text{ ${}^{\wedge}$
}d\text{y1}\text{ ${}^{\wedge}$
}d\text{z1}\text{ ${}^{\wedge}$ }d\text{z2}+(-1990656
\text{x1} \left(x^2+y^2\right)^6 \left(-2+x^4+2 x^2 y^2+y^4\right)^4\text{)
}d\text{t1}\text{ ${}^{\wedge}$ }dx\text{ ${}^{\wedge}$ }d\text{x1}\text{
${}^{\wedge}$ }d\text{x2}\text{ ${}^{\wedge}$ }dy\text{ ${}^{\wedge}$
}d\text{y2}\text{ ${}^{\wedge}$ }d\text{z1}\text{ ${}^{\wedge}$
}d\text{z2}+(-1990656
\left(x^2+y^2\right)^6 \left(-2+x^4+2 x^2 y^2+y^4\right)^4 \text{)
}d\text{t1}\text{
${}^{\wedge}$ }dx\text{ ${}^{\wedge}$ }d\text{x1}\text{ ${}^{\wedge}$
}d\text{x2}\text{ ${}^{\wedge}$ }dy\text{ ${}^{\wedge}$ }d\text{y1}\text{
${}^{\wedge}$
}d\text{y2}\text{ ${}^{\wedge}$ }d\text{z2}\),

\medskip

\NI As $\iota_Rdvol_{\mathbb{R}^9}=\tau^4$ for some $R \in
lin\{\frac{\partial}{\partial z_1},
\frac{\partial}{\partial z_2},\frac{\partial}{\partial
y_1},\frac{\partial}{\partial y_2}\},$
hence the Reeb field $R_{\tau}$ of $\tau$ is equal to $R$ because $\iota_R
\tau^4=
\iota_R\iota_Rdvol_{\mathbb{R}^9}=0.$ The field $R_{\tau}$ is obviously
perpendicular to
$\frac{\partial}{\partial r}$ (away from the degenerate set $\Sigma$).

\bibliographystyle{../../tex/macros/plain}
\bibliographystyle{amsalpha}

\medskip

\noindent {\bf BH:\ Mathematical Institute, Wroc\l aw University,

\noindent pl. Grunwaldzki 2/4,

\noindent 50-384 Wroc\l aw, Poland}

\NI and

\noindent {\bf  Department of Mathematics and Information
Technology,

\noindent University of Warmia and Mazury,

\noindent S\l oneczna 54, 10-710 Olsztyn, Poland}

\begin{flushleft}
\tt hajduk@math.uni.wroc.pl
\end{flushleft}
\medskip

\noindent {\bf RW: \ West Pomeranian University of Technology,

\noindent Mathematical Institute}

\noindent {\bf Al. Piast\'{o}w 48/49, 70--311 Szczecin, Poland}

\begin{flushleft}
 \tt rafal\_walczak2@wp.pl
\end{flushleft}

\end{document}